\documentclass{article}

\usepackage{amsmath, amsthm, amsfonts}

\usepackage[utf8]{inputenc}
\usepackage[T1]{fontenc}

\newtheorem{thm}{Theorem}

\newtheorem{prop}[thm]{Proposition}
\theoremstyle{definition}
\newtheorem{defn}[thm]{Definition}
\theoremstyle{remark}

\def\CC{\mathbb{C}}

\title{Limit holomorphic sections and Donaldson's  construction of symplectic submanifolds}
\author{Jean-Paul Mohsen \footnote{{\it jean-paul.mohsen@univ-amu.fr}}}

\begin{document}
\maketitle

Donaldson proved (in  \cite{Do96}) that if $L$ is a suitable positive line bundle over a closed symplectic manifold $X$,
then, for $k$ sufficiently large, the tensor power $L^k$ admits sections whose zero sets are symplectic submanifolds of $X$ 
(the sections are approximately holomorphic and they satisfy some uniform transversality condition).
The construction relies on the following observation: the local geometry of the bundles $L^k$ near any point $p\in X$,
after a normalization, converges to a model holomorphic Hermitian line bundle $K$ over (some ball in) the tangent space $T_p X$.
In this note, we will describe this phenomenon in detail and exploit it to reformulate Donaldson's theorem as a compactness result:
near each point $p$, the sections he obtains accumulate to holomorphic sections
of $K$ (that we call ``limit sections'')
and their uniform transversality  properties correspond to transversality properties of their limits.
Of course, similar considerations apply to all constructions based on Donaldson's techniques (e.g. \cite{Au01},  \cite{IbMaPr00}).

{\bf Acknowledgements.}
I want to thank Emmanuel Giroux for many important suggestions.

\section{Limit sections}

Let $X=(X,\,  \omega,\, J,\, g)$ be a closed almost-K\"{a}hler manifold. Hence $\omega$ is a symplectic form, $J$ is an almost-complex structure
and $g$ is a Riemannian metric, satisfying the following compatibility condition:
$g(V,W) = \omega(V,JW)$. 
Endow $X$ with a prequantization $L$ (a prequantization is a Hermitian line bundle over $X$ equipped
with a unitary connection of curvature $-i2\pi \omega$). 

The charts we will use are normal coordinates with respect to the renormalized metric $g_k=kg$.
Let $B\subset \mathbb{C}^n$ denote the unit ball, with  $n= \frac{1}{2}\, $dim$_{\mathbb{R}} X$.
Fix, for every large integer $k$, a chart $\varphi_k : B \rightarrow X$ satisfying two conditions:

(1) The chart $\varphi_k$ is an exponential map for the metric $g_k$ (i.e. given any unit vector $v\in \mathbb{C}^n$, the curve $t \mapsto \varphi_k (tv)$
is a geodesic with $g_k-$length $1$ velocity vector). 

(2) The differential $D\varphi_k (0)$ is a $\mathbb{C} - $linear map.
\\
\\   
Since $\varphi_k$ is a local diffeomorphism, one can transfer to $B$ the renormalized almost-K\"{a}hler structure 
$(\omega_k = k\omega,\, J,\, g_k =kg)$ and it is well known that 
this almost-K\"{a}hler structure tends to the standard flat K\"{a}hler structure on $B$, as $k \rightarrow \infty$, in the  
${\cal C}^{\infty} -$topology. 

The following observation is well-known to experts: the local geometry of
the bundle $L^k$ 
converges to the geometry of a model line bundle.
Fix some unitary radially flat isomorphism between the pullback line bundle $\varphi_k^* L^k$
and the trivial Hermitian line bundle $B \times \mathbb{C} \rightarrow B$.
Hence, the connection of $\varphi_k^* L^k$ induces a unitary connection $\nabla^k$ on $B \times \mathbb{C} \rightarrow B$.
As $k \rightarrow \infty$, the connection $\nabla^k$
tends to some model connection $\nabla^{\infty}$ on $B \times \mathbb{C} \rightarrow B$,
defined by:
$$
\nabla^{\infty} = d- i\pi \sum_{\alpha = 1}^{n} (x_{\alpha} dy_{\alpha} - y_{\alpha} dx_{\alpha}).
$$

There is a more conceptual description of $\nabla^{\infty}$: 
the model connection $\nabla^{\infty}$ is the only radially trivial connection with curvature 
$-i2 \pi \sum_{\alpha = 1}^{n} dx_{\alpha} \wedge dy_{\alpha}$.
\\
\\
Warning. Let $s$ be a section of the trivial bundle $B \times \mathbb{C} \rightarrow B$.
We say that $s$ is holomorphic if it is holomorphic for the connection $\nabla^{\infty}$.
Although the section $s$ is a function, it is not the usual concept of holomorphic function.
For example, the function $exp\left( - \frac{\pi}{2} \sum_{\alpha = 1}^{n}  | z_{\alpha} |^2 \right)$
is a holomorphic section and, more generally, the section $s$ 
of $B \times \mathbb{C} \rightarrow B$
is holomorphic if and only if the function 
$s \,exp\left( \frac{\pi}{2} \sum_{\alpha = 1}^{n}  | z_{\alpha} |^2 \right)$
is holomorphic in the usual sense.
\\
\\
This set of tools is well-known to experts. We will use it to study sequences of sections.
The following two definitions play an important role in our reformulation of Donaldson's theory.

\begin{defn}\label{Def1}
For every sufficiently large integer $k$, let $s_k$ be a ${\cal C}^{\infty} -$smooth section of $L^k$.
We say that the sequence $(s_k)$ is {\it renormalizable} if it satisfies the following compactness condition.

Let $(k_l)$ be a subsequence of the positive integers. 
For every sufficiently large integer
$l$, let $\varphi_l$ be a chart satisfying conditions 
(1) for $k_l$ (that is,
$\varphi_l$ is an exponential
map for $g_{k_l}$) 
and (2) and let $j_l$ be a
unitary radially flat isomorphism between the trivial line bundle  $B \times \mathbb{C} \rightarrow B$
and the pullback bundle $\varphi_l^* L^{k_l}$. 
If $\sigma_l$ denotes the section of the trivial
bundle  $B \times \mathbb{C} \rightarrow B$ corresponding to the pullback section $\varphi_l^*s_{k_l}$ via the
isomorphism $j_l$, 
then the sequence $(\sigma_l)$ has a subsequence $(\sigma_{l_m})$ which converges over $B$  for the  
smooth compact-open topology.
\end{defn}

\begin{defn}\label{Def2}
The limit of $(\sigma_{l_m})$ is called a {\it limit section} of the renormalizable sequence $(s_k)$.
Hence, a limit section is a section of  $B \times \mathbb{C} \rightarrow B$.
\end{defn}

We emphasize that we {\it don't} assume that all charts $\varphi_l$ have the same center.
\\
\\
Let $(s_k)$ be a renormalizable sequence.
If the sections $s_k$ are holomorphic then the limit sections are holomorphic.
More generally, let's state an informal principle:
if the sections $s_k$ satisfy some closed condition then one may
infer that the limit sections satisfy some corresponding condition.
We won't be more specific about this principle (we won't even explain the meaning of the word {\it closed}).

Of course, concerning open conditions, it goes in the opposite direction.
For example, if all limit sections are transverse to $0$ then $s_k$ is transverse to $0$, for every sufficiently large integer $k$.
If, in addition, the zero sets of limit sections are symplectic submanifolds in $B$
(for the symplectic form $\sum_{\alpha = 1}^{n} dx_{\alpha} \wedge dy_{\alpha}$),
then for every sufficiently large integer $k$, the zero set of $s_k$ is a symplectic submanifold in $X$.
Note that every complex submanifold is symplectic. Hence one get the following proposition.

\begin{prop}{\label{o}}
For every sufficiently large integer $k$, let $s_k$ be a ${\cal C}^{\infty} -$smooth section of $L^k$.
Suppose that $(s_k)$ is a renormalizable sequence and suppose that every limit section of 
the sequence $(s_k)$ is holomorphic and transverse to $0$.
Then for every sufficiently large integer $k$, the zero set of $s_k$ is a codimension $2$ symplectic submanifold in $X$.
\end{prop}

In the integrable case ($X$ is K\"{a}hler), the sections $s_k$ we will consider
are often holomorphic whereas in the non-integrable case ($X$ almost-K\"{a}hler), typically,
the limit sections are holomorphic but the sections $s_k$ aren't. 

\section{Transversality theorems}

To compare with, let's state a consequence of the Kodaira embedding theorem.

\begin{thm}{\label{BK}}
Suppose $J$ is integrable. Then, for every sufficiently large integer $k$, there exists a holomorphic section $s_k$ of $L^k$ which is transverse to $0$.
\end{thm}
Proof.
Kodaira's theorem implies that, for every sufficiently large $k$, there are no base points. 
Hence, almost every section is transverse to $0$, by the Bertini theorem. 
\\
\\
\indent
Let's first state Donaldson's theorem in the integrable case.

\begin{thm}{\label{DI}}
Suppose $J$ is integrable. Then, for every integer $k\geq 1$, there exists a holomorphic section $s_k$ of $L^k$ such that:

(1) The sequence $(s_k)$ is renormalizable.

(2) The limit sections of the sequence $(s_k)$ are transverse to $0$.
\end{thm}
In the integrable case, one may describe Donaldson's theorem an elaborate variant of Theorem {\ref{BK}}.
The variant has the advantage of being easily transferable to symplectic geometry.
Of course this was Donaldson's main goal and most applications of his techniques are symplectic and contact results. 
It is known that if $X$ is almost-K\"{a}hler, then, in general, one can't get holomorphic sections.
Nevertheless, we get {\it asymptotically holomorphic} sections.
In our reformulation, the definition is quite simple:
a renormalizable sequence of smooth sections
is {\it asymptotically holomorphic} if every limit section is holomorphic for the connection $\nabla^{\infty}$.
Note that this definition is weaker than the usual quantitative definition. However, in practice, 
the following version of Donaldson's theorem is sufficient for many corollaries.

\begin{thm}{\label{D}}
For every integer $k\geq 1$, there exists a ${\cal C}^{\infty} -$smooth section $s_k$ of $L^k$ such that:

(1) The sequence $(s_k)$ is renormalizable.

(2) The limit sections of the sequence $(s_k)$ are holomorphic and transverse to $0$.
\end{thm}

(Hence, for every sufficiently large integer $k$, the section $s_k$ is transverse to $0$ and, by Proposition \ref{o},
the zero set of $s_k$ is a codimension $2$ symplectic submanifold.)
\\
\\
\indent
Proof of Theorem {\ref{DI}} and Theorem {\ref{D}}.
Donaldson's techniques (in  \cite{Do96}, see also \cite{Au97} and \cite{Do99}) produce sections $s_k$ satisfying two famillies of estimates: 
\begin{eqnarray*}
\| s_k \|_{{\cal C}^r, g_k} & = & O (1)
\\
\| \overline{\partial} s_k \|_{{\cal C}^r , g_k} & = & O (k^{-\frac{1}{2}})
\end{eqnarray*}
(for every natural integer $r$),
and a uniform transversality condition:
\begin{eqnarray*}
\min_{p\in X} \left( \| s_k(p) \| + \| \nabla s_k(p) \|_{g_k} \right) & \geq & \eta
\end{eqnarray*}
where one calculates the ${\cal C}^r-$norm and the norm of $\nabla s_k(p)$  with the renormalized metric $g_k = kg$.
Here $\eta$ denotes a positive number, independent of $k$. 

Recall the notations of Definition {\ref{Def1}}.
The section $\sigma_l$ of  $B \times \mathbb{C} \rightarrow B$
corresponds to the pull-back section $\varphi_l^*s_{k_l}$ where $\varphi_l$ is the exponential map for the renormalized metric $g_{k_l}$.
The first estimate implies the following estimate:
$$
\| \sigma_{l} \|_{{\cal C}^r}  =  O (1)
$$
on the unit ball $B$. Hence, some subsequence of $(\sigma_l)$ converges in the smooth topology and the sequence $(s_k)$ is renormalizable.

The connection of $\varphi_l^* L^{k_l}$ induces a unitary connection $\nabla^{k_l}$ on $B \times \mathbb{C} \rightarrow B$.
Let  $\overline{\partial}^{k_l}$ denote the $(0,1)-$part of $\nabla^{k_l}$
and let $\overline{\partial}^{\infty}$ denote the $(0,1)-$part of the limit connection $\nabla^{\infty}$.
Donaldson's second estimate implies the following estimate:
$$
\|  \overline{\partial}^{k_l} \sigma_{l} \|_{{\cal C}^r}  =  O (k_l^{-\frac{1}{2}}).
$$
Since $\nabla^{k_l}$ tends to  $\nabla^{\infty}$, the $(0,1)-$part
$\overline{\partial}^{\infty} \sigma_l$ tends to $0$ and the sequence $(s_k)$ is asymptotically holomorphic.

The third estimate implies the following estimate:
$$
\min_{p\in B} \left( \| \sigma_{l}(p) \| + \| \nabla \sigma_{l}(p) \| \right) \geq  \frac{\eta}{2}.
$$
Every limit $\sigma_{\infty}$ of a subsequence $(\sigma_{l_m})$ (see Definition \ref{Def2}) satisfies the same estimate:
$$
\min_{p\in B} \left( \| \sigma_{\infty}(p) \| + \| \nabla \sigma_{\infty}(p) \| \right)  \geq  \frac{\eta}{2}.
$$
Hence, $\sigma_{\infty}$ is transverse to $0$.
The proof of Theorem \ref{D} is completed. 
In the integrable case, Donaldson's sections are holomorphic and the proof of Theorem \ref{DI} is similar.
\\
\\
\indent
As noted by Donaldson, the asymptotic transversality property provides bounds for the Riemannian geometry of the zero set, see
\cite[Corollary 33]{Do96}.
For example, one gets the following result.

\begin{prop}
Let $(s_k)$ be a renormalizable sequence. Suppose every limit sequence is transverse to $0$.
Let $Y_k$ be the zero set of $s_k$.
For every sufficiently large integer $k$,
if $p$ lies in $Y_k$ and $A\subset T_p Y_k$ is a 2-plane, 
then the sectional curvature $K_{Y_k,g_k} (p,A)$ of $Y_k$ at $(p,A)$ for the metric $g_k$
satisfies the following estimate:
$$
 | K_{Y_k,g_k} (p,A) | \leq C
$$
where the bound $C$ is independent of $k$, $p$ and $A$.
\end{prop}

(Hence, if one prefers to calculate with the metric $g$, the sectional curvature is
bounded by some linear function of $k$ because $K_{Y_k,g} (p,A) = k\, K_{Y_k,g_k} (p,A)$.)
\\
\\
\indent
Proof.
Define $u_k= \max_{(p,A)} | K_{Y_k,g_k} (p,A) |$. Since $Y_k$ is compact,
there exist a point $p_k \in Y_k$ and a 2-plane $A_k \subset T_{p_k} Y_k$ satisfying the following equation:
$$ | K_{Y_k,g_k} (p_k,A_k ) | = u_k.$$

Let $\varphi_k$ be a chart centered at $p_k$ satifying conditions (1) and (2)
of Definition {\ref{Def1}}. 
Consider the 2-plane $A_k' =(d\varphi_k(0))^{-1} (A_k) \subset \CC^n$.
 Since the set of 2-planes of $\CC^n$ is compact, every subsequence of $(A_k')$
 admits a convergent subsubsequence $(A_{k_l}')$. 
 The corresponding sequence $\sigma_l$ 
 (using the notations of Definition {\ref{Def1}})
 admits a limit $\sigma_{\infty}$ which is transverse to $0$. Therefore 
the zero set $Y_{\infty} \subset B$ of $\sigma_{\infty}$
is a submanifold and the local geometry of
the corresponding submanifolds $(Y_{k_l})$
converges to the geometry of $Y_{\infty}$.
In particular, the sequence
$(u_{k_l})$ tends to $| K_{Y_{\infty},\mu} (0,A'_{\infty}) |$
where $\mu$ is the standard Euclidean metric on $\CC^n$
and $A'_{\infty}$ is the limit 2-plane.

Hence every subsequence of the sequence $(u_k)$ admits a convergent subsubsequence and therefore $(u_k)$ is a bounded sequence.




\end{document}